%
%
\input amstex
\documentstyle{amsppt}
\TagsOnRight
\magnification=\magstep1
%
%
\parindent=12pt                   
\parskip=0pt
%
%
%
%
\loadbold
\loadeusm
\def\scr#1{{\fam\eusmfam\relax#1}}
\NoBlackBoxes
%
\def\lines#1:{\ifhmode\vskip#1\baselineskip\noindent\leftskip=0pt
              \else\ifvmode\vskip#1\baselineskip
                    \fi
              \fi}
\def\openC{\Bbb C}

\def\openH{\Bbb H}
\def\openN{\text{I\!N}}
\def\openF{\text{I\!F}}
\def\openP{\Bbb P}
\def\Pn#1{\openP^#1}
\def\kdual{\scriptscriptstyle{\vee}}
\def\Pnd#1{{\check\openP}^#1}

\def\openZ{\Bbb Z}

\def\coker{\operatorname{coker}\,}
\def\rk{\operatorname{rk}\,}

\def\ker{\operatorname{ker}\,}
\def\id{\operatorname{id}\,}
\def\Hom{\operatorname{Hom}}
\def\Ext{\operatorname{Ext}}
\def\deg{\operatorname{deg}\,}

\def\Dach{\operatorname{\Lambda}\,}

\def\dsum#1#2{\operatornamewithlimits{\oplus}_{#1}^{#2}}
\def\qed{}
%
%
\def\id#1#2{{\scr I}_{#1}(#2)}
%
\def\struct#1#2{{\scr O}_{#1}(#2)}
%

%
\def\systruct#1{{\scr O}(#1)}
%
\def\somega#1#2{\Omega^{#1}(#2)}
%
\def\mapleft#1{\overset{#1}\to\longleftarrow}
\def\mapright#1{\overset{#1}\to\longrightarrow}
\def\eso#1#2#3#4#5#6{0\longrightarrow \; {\scr O}_{#1}(#2)\; \longrightarrow \;
{\scr O}_{#3}(#4)\;
\longrightarrow \; {\scr O}_{#5}(#6) \; \longrightarrow 0}
%
\def\esi#1#2#3#4#5#6{0\longrightarrow \;  {\scr I}_{#1}(#2)\; \longrightarrow
\;
{\scr I}_{#3}(#4)\; \longrightarrow \; {\scr I}_{#5}(#6) \; \longrightarrow 0}
%

%

%

%
\def\exactg#1#2#3{0\longrightarrow \;{#1}\; \longrightarrow \;
{#2}\; \longrightarrow \;  {#3}\; \longrightarrow 0}
%
\def\sexactg#1#2#3#4#5{0\longrightarrow \;{#1}\; \mapright{#2} \;
{#3}\; \mapright{#4} \;  {#5}\; \longrightarrow 0}
%
\def\hi#1{{\text h}^{#1}}
\def\Hi#1{{\text H}^{#1}}
%

%
%

%
\def\sCohid#1#2#3{{\text H}^{#1}({\scr I}_{#2}({#3}))}
%
%
\def\sCohos#1#2#3{{\text H}^{#1}({\scr O}_{#2}({#3}))}
%

%
%

%
\def\scohid#1#2#3{{\text h}^{#1}({\scr I}_{#2}({#3}))}
%
\def\scohos#1#2#3{{\text h}^{#1}({\scr O}_{#2}({#3}))}
%

%
\def\sCoh#1#2#3{{\text H}^{#1}({#2}({#3}))}
%

%
%
\topmatter
\title
Examples of smooth non-general type surfaces in ${\openP}^4$
\endtitle
\author
Sorin Popescu
\endauthor
%
%
\address
Columbia University,
Department of Mathematics,
Mail Code 4417,
2990 Broadway,
New York, NY 10027,
USA
\endaddress
\email
psorin\@math.columbia.edu
\endemail
\thanks
Much of this paper is part of the author's  Ph. D. thesis. The work on it
was initiated during stays at the University of Bayreuth and the University
of Kaiserslautern, to whom also I thank for financial support. Special thanks
go to my advisors W. Decker and F.-O. Schreyer whose support and
advice guided me all along the elaboration of my thesis.
Also I would like to thank A. Aure, W. Barth, G. Ellingsrud,
Ch. Peskine and K. Ranestad for many
inspiring conversations. The author has lectured on this paper
at the Europroj Conference, Barcelona 1994.
\endthanks

%
%
\endtopmatter
%

\long\def\beiliw#1#2#3#4#5#6#7#8{
\vbox{\tabskip=0pt\offinterlineskip
\def\tablerule{\omit&\multispan{11}{\hrulefill}&&\cr}
\halign to250pt{\strut##  & \vrule##\tabskip=1em  &\hbox to
0.19cm{\hfil##\hfil}&\vrule##&
              \hbox to 0.19cm{\hfil##\hfil}& \vrule##&\hbox to
0.19cm{\hfil##\hfil}&\vrule##&
              \hbox to 0.19cm{\hfil##\hfil}& \vrule##&\hbox to
0.19cm{\hfil##\hfil}&\vrule##&
              \hbox to 0.19cm{\hfil##\hfil}& ##\hfil &\hfil## \tabskip=0em \cr
\hskip-13pt
i&\omit\hskip-2.2pt$\uparrow$&&\omit&&\omit&&\omit&&\omit&&\omit&&&\cr
\noalign{\vskip-2pt}
&\omit\vrule height10pt&\omit &\omit&&\omit&&\omit&&\omit&&\omit&&&\cr%
\tablerule
\omit&height5pt&\omit&&\omit&&\omit&&\omit&&\omit&&\omit&\omit&\cr
&& && && && && &&&&\cr
\tablerule
\omit&height5pt&\omit&&\omit&&\omit&&\omit&&\omit&&\omit&\omit&\cr
&&#1&& && && && &&&&\cr
\tablerule
\omit&height5pt&\omit&&\omit&&\omit&&\omit&&\omit&&\omit&\omit&\cr
&&#2&&#3&&#4&&#5&& &&&\quad\quad\quad$\scohid i S p$&\cr
\tablerule
\omit&height5pt&\omit&&\omit&&\omit&&\omit&&\omit&&\omit&\omit&\cr
&& && &&#6&&#7&& &&&&\cr
\tablerule
\omit&height5pt&\omit&&\omit&&\omit&&\omit&&\omit&&\omit&\omit&\cr
&& && && && &&#8&&&&\cr
\noalign{\vskip-3.1pt}
\omit&\multispan{12}{\rightarrowfill}&\cr
\noalign{\vskip5pt}
\multispan{13}{\hfil p}&\cr
}
}
}
\long\def\beiliv#1#2#3#4#5#6#7#8{
\vbox{\tabskip=0pt\offinterlineskip
\def\tablerule{\omit&\multispan{11}{\hrulefill}&&\cr}
\halign to250pt{\strut##  & \vrule##\tabskip=1em  &\hbox to
0.19cm{\hfil##\hfil}&\vrule##&
              \hbox to 0.19cm{\hfil##\hfil}& \vrule##&\hbox to
0.19cm{\hfil##\hfil}&\vrule##&
              \hbox to 0.19cm{\hfil##\hfil}& \vrule##&\hbox to
0.19cm{\hfil##\hfil}&\vrule##&
              \hbox to 0.19cm{\hfil##\hfil}& ##\hfil &\hfil## \tabskip=0em \cr
\hskip-13pt
i&\omit\hskip-2.2pt$\uparrow$&&\omit&&\omit&&\omit&&\omit&&\omit&&&\cr
\noalign{\vskip-2pt}
&\omit\vrule height10pt&\omit &\omit&&\omit&&\omit&&\omit&&\omit&&&\cr%
\tablerule
\omit&height5pt&\omit&&\omit&&\omit&&\omit&&\omit&&\omit&\omit&\cr
&& && && && && &&&&\cr
\tablerule
\omit&height5pt&\omit&&\omit&&\omit&&\omit&&\omit&&\omit&\omit&\cr
&&#1&& && && && &&&&\cr
\tablerule
\omit&height5pt&\omit&&\omit&&\omit&&\omit&&\omit&&\omit&\omit&\cr
&&#2&&#3&&#4&& && &&&\quad\quad\quad$\scohid i S p$&\cr
\tablerule
\omit&height5pt&\omit&&\omit&&\omit&&\omit&&\omit&&\omit&\omit&\cr
&& && &&#5&&#6&&#7&&&&\cr
\tablerule
\omit&height5pt&\omit&&\omit&&\omit&&\omit&&\omit&&\omit&\omit&\cr
&& && && && &&#8&&&&\cr
\noalign{\vskip-3.1pt}
\omit&\multispan{12}{\rightarrowfill}&\cr
\noalign{\vskip5pt}
\multispan{13}{\hfil p}&\cr
}
}
}

\long\def\beiliz#1#2#3#4#5#6#7#8{
\vbox{\tabskip=0pt\offinterlineskip
\def\tablerule{\omit&\multispan{11}{\hrulefill}&&\cr}
\halign to250pt{\strut##  & \vrule##\tabskip=1em  &\hbox to
0.19cm{\hfil##\hfil}&\vrule##&
              \hbox to 0.19cm{\hfil##\hfil}& \vrule##&\hbox to
0.19cm{\hfil##\hfil}&\vrule##&
              \hbox to 0.19cm{\hfil##\hfil}& \vrule##&\hbox to
0.19cm{\hfil##\hfil}&\vrule##&
              \hbox to 0.19cm{\hfil##\hfil}& ##\hfil &\hfil## \tabskip=0em \cr
\hskip-13pt
i&\omit\hskip-2.2pt$\uparrow$&&\omit&&\omit&&\omit&&\omit&&\omit&&&\cr
\noalign{\vskip-2pt}
&\omit\vrule height10pt&\omit &\omit&&\omit&&\omit&&\omit&&\omit&&&\cr%
\tablerule
\omit&height5pt&\omit&&\omit&&\omit&&\omit&&\omit&&\omit&\omit&\cr
&& && && && && &&&&\cr
\tablerule
\omit&height5pt&\omit&&\omit&&\omit&&\omit&&\omit&&\omit&\omit&\cr
&&#1&& && && && &&&&\cr
\tablerule
\omit&height5pt&\omit&&\omit&&\omit&&\omit&&\omit&&\omit&\omit&\cr
&& &&#2&&#3&&#4&& &&&\quad\quad\quad$\scohid i S p$&\cr
\tablerule
\omit&height5pt&\omit&&\omit&&\omit&&\omit&&\omit&&\omit&\omit&\cr
&& && &&#5&&#6&&#7&&&&\cr
\tablerule
\omit&height5pt&\omit&&\omit&&\omit&&\omit&&\omit&&\omit&\omit&\cr
&& && && && &&#8&&&&\cr
\noalign{\vskip-3.1pt}
\omit&\multispan{12}{\rightarrowfill}&\cr
\noalign{\vskip5pt}
\multispan{13}{\hfil p}&\cr
}
}
}
\parindent=0pt
\document
\head{0. Introduction}\endhead
Smooth projective  varieties with small invariants have got renewed
interest in recent years,
primarily due to the fine study of the adjunction mapping by Reider, Sommese,
Van de Ven and others. For the special case of smooth surfaces in $\Pn 4$
the method goes back to the Italian geometers, who at the turn of the
century used it
for the study of the surfaces of degree less than $7$,
or sectional genus $\pi\le 3$.
Later on, for larger values of the  invariants,
there are contributions by Commesatti and especially
Roth. For example, in \cite{38}, Roth tried to establish a
classification of smooth surfaces with
$\pi\le 6$, but his lists are incomplete since
he misses the non-special rational
surfaces of degree $9$ and the minimal bielliptic
surfaces of degree $10$. Nowadays,
through the effort of several mathematicians
(some references are given below), a complete
classification of smooth surfaces in $\Pn 4$ has been worked out up to
degree $10$, and a partial one is available in degree $11$.
\par\medskip

But, apart from the general framework of classification problems
concerning codimension
two varieties, there is another strong motivation for the interest in
these surfaces. Namely, in a recent paper Ellingsrud and Peskine
\cite{19} proved  Hartshorne's conjecture
that there are only finitely many families of special surfaces in
$\Pn 4$. More specifically,
given an integer $a<6$, they show that the degree of smooth
surfaces with  $K^2\le a\chi$ is
bounded. In particular, there are only finitely many
families of smooth surfaces in $\Pn 4$,
not of general type. However, the question of an exact degree
bound is still open. A recent
work of Braun and Fl\o ystad \cite{10} improves the initial
bound ($\sim 10000$) of Ellingsrud and Peskine
to $d\le 105$, but it is believed that the
degree of the smooth, non-general type surfaces
in $\Pn 4$ should be less than or equal to $15$. A similar
finiteness result for $3$-folds in
$\Pn 5$ was proved in \cite{11}, but the real degree bound is believed
to be much higher
in this case. Nevertheless, examples of smooth 3-folds in $\Pn 5$ not of
general type  are known only up to degree $18$ (see \cite{18}
for more details and a complete list of known examples).\par\medskip

Another reason for the interest in studying surfaces in $\Pn 4$ is
the small number  of known liaison classes of such surfaces.
Each new specimen of liaison classes is of  real
interest. In this direction, the work of
Decker, Ein and Schreyer \cite{17} provides a powerful
and effective method of construction of surfaces in $\Pn 4$.
\par\medskip

The aim of this paper is to provide a series of examples of
smooth surfaces in $\Pn 4$,
not of general type, in degrees varying from $12$ up to $14$,
and to describe part of their geometry.
In degree 15, two families of abelian surfaces \cite{23}, \cite{3}, \cite{34},
and a family of bielliptic surfaces \cite{4} are currently known.
We have tried to work out examples of degrees higher then 15 but
failed  in this attempt.
The methods of construction we used are mainly the syzygy approach of \cite{17}
and liaison techniques. The families we construct are:\par\medskip

\item\item{-} minimal proper elliptic surfaces of degree
$12$ and sectional genus $\pi=13$,\smallskip
\item\item{-} two types of non-minimal proper elliptic surfaces of degree $12$
and sectional genus $\pi=14$,\smallskip
\item\item{-} non-minimal $K3$ surfaces of degree $13$ and sectional genus
$16$, and\smallskip
\item\item{-} non-minimal $K3$ surfaces of degree $14$ and sectional genus
$19$.
\par\medskip

At this point it may be appropriate to recall some references for the list of
the smooth surfaces in $\Pn 4$ of degrees less or equal to $11$.
The classification and construction of surfaces of degree
$\le 7$ was initiated in \cite{38} and completed up to degree $8$ in \cite{24},
\cite{25}, \cite{30}, \cite{31}, \cite{32},
supplemented for the case of rational surfaces of degree $8$, sectional genus
$5$ by \cite{1}. In
degree $9$, the rational surfaces are described in \cite{1} and \cite{2},
the Enriques surfaces with $\pi=6$ in \cite{14} and \cite{15}, while the
classification and description of the liaison classes is
completed in \cite{5}. In degree $10$, the classification in terms
of numerical invariants and the description of a large number of
surfaces is achieved in the beautiful
thesis of K. Ranestad \cite{36}. The existence, the uniqueness and
the geometry of bielliptic surfaces
of degree $10$, $\pi=6$ were taken care by \cite{42} and \cite{4},
respectively,
the Enriques surfaces of degree $10$, $\pi=8$ were first constructed in
\cite{17} and further studied in \cite{12}, while the minimal
abelian surfaces were first described about
$20$ years ago in \cite{23}.  Finally, a non-minimal $K3$
surface of degree $10$,
$\pi=9$, lying on only one quartic hypersurface was
constructed \cite{34},
thus giving a positive answer to the last open existence case in
\cite{36}.
The remaining uniqueness problems, the syzygies and the description
of the liaison classes in
degree $10$ were completed in \cite{35}. Finally, 23 different
families of smooth surfaces
are known in degree 11, and \cite{34} is an attempt at construction
and classification  in this degree.
\par\bigskip

\head {1. Preliminaries}\endhead
\subhead Double point formula\endsubhead
For smooth surfaces $S\subset\Pn 4$  the relation \cite{22}:
$$d^2-c_2(N_S)=d^2-10d-5H\cdot K-2K^2+12\chi(S)=0,\leqno(1.1)$$
expresses the fact that $S$ has no double points.
\subhead Linkage \cite{33}\endsubhead
Two surfaces $S$ and $S^{\prime}$ in $\Pn 4$ are
said to be linked $(m,n)$ if there exist hypersurfaces
$V$ and $V^{\prime}$ of
degree $n$ and $m$ respectively, without
common components and  such that $V\cap V^{\prime}=S\cup S^{\prime}$.
The standard sequences of linkage, namely
$$\eso SK{S\cup S^{\prime}}{m+n-5}{S^{\prime}}{m+n-5}$$
$$\eso SKS{m+n-5}{S\cap S^{\prime}}{m+n-5},$$  yield then
$\chi (S^{\prime})=\chi (V\cap V^{\prime})-\chi (\struct S{m+n-5})$ and
a relation between the sectional genera:
$\pi (S)-\pi (S^{\prime})={1\over 2}(m+n-4)(d(S)-d(S^{\prime})).$

\subhead The Eagon-Northcott complex method \cite{17}\endsubhead
The aim of this method is to realize a surface $S\subset\Pn 4$ as the
determinantal locus
$S = D(\varphi) = \lbrace \; x\in \Pn 4\; \mid\; {\rk}\; \varphi(x)<e\;\rbrace$
of a map $\varphi$ between two vector bundles
${\scr E}$ of rank $e$, and ${\scr F}$ of rank $e+1$ on $\Pn 4$.
In  case $\varphi$ degenerates
in (expected) codimension two, $D(\varphi)$ is locally Cohen-Macaulay
and the Eagon-Northcott complex is exact and identifies
$\coker\varphi$ with the twisted ideal sheaf of $S$:
$$\sexactg {\scr E} {\varphi} {\scr F} {}
{\id S {c_1({\scr F})- c_1({\scr E})}}.$$
In order to construct a surface $S$ with the desired
invariants one has to find appropriate
vector bundles $\scr E$ and $\scr F$, and a general method for this
is to determine first
the differentials of the Beilinson spectral sequence applied to
the twisted ideal sheaf $\id S m$,
for a suitable $m\in\openN$:
\proclaim {Theorem 1.1.\cite{9}} Let ${\scr G}$ be a coherent sheaf on $\Pn n =
\openP(V)$. There
exists a spectral sequence with $E_1$ terms
$$E_1^{pq} = {\text H}^q(\Pn n , {\scr G}(p))\otimes {\Omega}_{\Pn
n}^{-p}(-p)$$
converging to ${\scr G}$; i.e. $E_{\infty}^{pq} = 0$ for $p + q \ne 0$ and
$\oplus E_{\infty}^{p,-p}$ is the associated graded sheaf of a suitable
filtration of
${\scr G}$.\endproclaim
All the $E_1$-terms are in the second quadrant and only finitely many of them
are non-zero. Via canonical isomorphisms induced by contraction,
$\Hom(\Omega_{\Pn n}^{i}(i),\Omega_{\Pn n}^{j}(j))\cong \Lambda^{i-j} V$, for
$i\ge j$,
cf. \cite{9},
the $d_1$-differentials
$$d_1^{pq}\in \Hom\Big({\text H}^q(\Pn n , {\scr G}(p))\otimes
{\Omega}_{\Pn n}^{-p}(-p),
{\text H}^q(\Pn n , {\scr G}(p+1))\otimes {\Omega}_{\Pn n}^{-p-1}(-p-1)\Big)$$
can be identified with the natural multiplication maps in
$$\Hom\Big({\text H}^0(\Pn n ,
{\scr O}_{\Pn n}(1))\otimes{\text H}^q(\Pn n , {\scr G}(p)),
{\text H}^q(\Pn n , {\scr G}(p+1))\Big).$$
In our specific case, this  means that to determine the $d_1$-differentials is
equivalent to fixing the
module structure of the Hartshorne-Rao modules $\oplus_{p}{\text H}^q
(\Pn 4 , {\id S p})$, $q\in\{1,2\}$.

\subhead Multisecants \cite{28}\endsubhead
 Some classical numerical formulas for multisecant lines to a smooth
surface $S$ in $\Pn 4$ have been recently given a modern treatment
by Le Barz. Consider the
double curve $\Gamma$ of a general projection of such a surface $S$ to $\Pn 3$
and denote by
$$\delta = {{d-1}\choose 2} - \pi$$
the degree of $\Gamma$, by
$$t = {{d-1}\choose 3} - \pi (d - 3) + 2\chi - 2$$ the number of apparent
triple points,
i.e., the number of trisecants to $S$ which meet a general point, and by
$$h = {1 \over 2}(\delta (\delta - d + 2) -3t)$$ the number of apparent double
points on $\Gamma$.
Suppose there are no lines on $S$ with positive self-intersection. Then the
number of 6-secants (if finite)
plus the number of exceptional lines on $S$ is:
$$\leqalignno{ N_6(d,\pi,\chi) =& -{1\over 144}d(d - 4)(d - 5)(d^3 + 30d^2 -
577d + 786) &(1.2)\cr
& + \delta(2{d\choose 4} + 2{d\choose 3} - 45{d\choose 2} + 148d - 317) &\cr
& - {1\over 2}{\delta\choose 2}(d^2 - 27d + 120) - 2{\delta\choose 3} &\cr
& + h(\delta - 8d + 56) + t(9d - 3\delta - 28) + {t\choose 2}. &\cr}$$

\head{2. Three families of smooth elliptic surfaces of degree 12
and a smooth $K3$ surface of degree 13}\endhead
We construct in the sequel several examples of smooth, regular, proper elliptic
surfaces
of degree $12$ in $\Pn 4$. Previously known examples of non-general type
surfaces with $d=12$ were only the blown-up $K3$ surfaces with $\pi=14$, with
one
exceptional quartic and ten exceptional lines constructed in
\cite{17}.\par\medskip
Recall first that for a smooth surface $S$ of degree $12$ in $\Pn 4$ the double
point formula reads
$$K^2=47-5\pi+6\chi,$$ while Severi's theorem \cite{43} and Riemann-Roch give
$$\pi=\chi+8+\scohos 1SH -\scohos 0S{K-H}.$$\bigskip

\proclaim {Proposition 2.1} There exist smooth, regular, minimal proper
elliptic surfaces
$S\subset\Pn 4$, with $d=12$, $\pi=13$, $\chi=3$, and  $10$ skew $6$-secant
lines.
\endproclaim
\demo{Proof} For construction we use the syzygy approach in \cite{17}. A
promising Beilinson
cohomology table is
$$\qquad\beiliv 2 \ 2 \ \ 4 5 \ $$
where ${\text h}^1({\scr I}_S)={\text h}^1({\scr I}_S(1))=0$ ($S$ being reduced
and linearly
normal, cf \cite{43}), ${\text h}^2({\scr I}_S)={\text h}^1({\scr O}_S)=q=0$,
${\text h}^3({\scr I}_S)={\text h}^2({\scr O}_S)=p_g=2$, ${\text h}^3({\scr
I}_S(k))=0$, for
$k\ge 1$ (since $\kappa(S)<2$), while ${\text h}^0({\scr I}_S(m))=0$ for $m\le
3$.
We may set ${\scr E}:=2{\scr O}(-1)\oplus 2\somega 3 3$ and ${\scr F}:
=\widetilde{(\ker\psi)}$,
where $4R(1)\mapleft{\psi}15R$ is the minimal free presentation of the graded
$R=\openC[x_0,x_1,x_2,x_3,x_4]$-module $\sCohid 1 S {\ast+4}
= \sCoh 1 {\scr F} {\ast}$,
and try to construct the surface as the degeneracy locus of a general morphism
$\varphi\in\Hom({\scr E},{\scr F})$.
However, for a general choice of the matrix $\psi$, the module
$M:=\coker\psi$ has a minimal free resolution of type
{\eightpoint
$$
\vbox{
\halign{&\hfil$\,#\,$\hfil\cr
M&\leftarrow&4R(1)&\mapleft{\psi}&15R&&15R(-1)\cr
&&&&&\vbox  to 10pt{\vskip-4pt\hbox{\hss$\nwarrow$}\vss}&\oplus\cr
&&&&&&10R(-2)&\vbox  to
10pt{\vskip-4pt\hbox{\hss$\nwarrow$}\vss}&30R(-3)&\leftarrow&21R(-4)&
\leftarrow&5R(-5)&\leftarrow0\cr
}}$$}
and thus $\Hom(\somega 3 3,{\scr F})=0$ in this case.  What is needed for the
construction to work
is that  $\psi$ has at least two linear syzygies of second order. We will
choose
a $\psi$ featuring such syzygies.\par\medskip

Let $F$ be the Horrocks-Mumford bundle (see \cite{23}). It is a stable rank $2$
vector bundle on $\Pn 4$ with Chern classes $c_1=-1$, $c_2=4$, and its ${\text
H}^1$-cohomology
module has a minimal free resolution of type (cf. \cite{16}, see also
\cite{34})
{\eightpoint
$$
\vbox{\vfil
\halign{&\hfil$\,#\,$\hfil\cr
0\leftarrow&\sCoh 1 F {\ast}&\leftarrow&5R&\mapleft{\gamma}&15R(-1)&&10R(-2)\cr
&&&&&&\vbox  to 10pt{\vskip-4pt\hbox{\hss$\nwarrow$}\vss}&\oplus\cr
&&&&&&&4R(-3)&&2R(-3)\cr
&&&&&&&\oplus&\vbox  to 10pt{\vskip-4pt\hbox{\hss$\nwarrow$}\vss}&\oplus\cr
&&&&&&&15R(-4)&&35R(-5)&\leftarrow&20R(-6)\cr
&&&&&&&&&&&&\vbox  to
10pt{\vskip-4pt\hbox{\hss$\nwarrow$}\vss}&2R(-8)&\leftarrow&0,\cr
}}
$$}
with $\gamma=\bigl( M_{z_0}(x)\mid M_{z_1}(x)\mid M_{z_2}(x)\bigr)$, where
$M_z(x)={(x_{i+j}z_{i-j})}_{{i,j}\in\openZ_5}$ are Moore blocks, and
the parameters are say $z_0=(1:0:0:0:0)$, $z_1=(0:1:0:0:1)$ and
$z_2=(0:0:1:1:0)$.\par\smallskip

{}From \cite{23}, or even just by looking at the above resolution
since $F(-1)$ is the cokernel
of the unique morphism $0\rightarrow 2\somega 3 2
\mapright{\theta}{Syz}_1(\sCoh 1 F {\ast})$, we have
$\hi 0 (F)=\hi 0 ({F(1)})=0$ and $\hi 0 {(F(2))}=4$.\par\smallskip

Consider now a rank $3$ vector bundle $E$ on $\Pn 4$ constructed
as the extension
$$\exactg F E {\scr O}\leqno{(2.2)}$$
corresponding to a non-trivial element
$0\ne\xi\in\Hi 1 (F)=\Ext^1({\scr O},F)$. Then $E$ has Chern
classes $c_1(E)=-1$, $c_2(E)=4$ and $c_3(E)=0$, and is stable because $F$ is.
Also $\hi 1(E)=
\hi 1(F)-1=4$ by construction. We will assume in the sequel that the extension
$E$ comes from
a generic element $\xi\in\Hi 1 (F)$, meaning by this that $\xi$ satisfies the
following
two conditions:\par\medskip
- the natural map $\openC\,\xi\otimes\sCohos 0
{\Pn 4} 1\rightarrow\sCoh 1 F 1$ induced by multiplication with
linear forms is injective, while\hfill\break
- the similar natural map $\openC\,\xi\otimes\sCohos 0
{\Pn 4} 2\rightarrow\sCoh 1 F 2$ is surjective.
\par\medskip

To see that such a choice is possible one either checks
it for a random $\xi$ via \cite{7}, or one uses
the invariance of $\sCoh 1F\ast$ under the group $G={\openH}_5\rtimes\openZ_2$,
where ${\openH}_5$
is the Heisenberg group of level $5$. Namely, using notations and
facts from \cite{23}, or \cite{29}, if
$\Pn 4=\openP(V)$ then we can find a basis
$e_0,\dots ,e_4$ of $V$ such that, under the Schr\"odinger
representation of ${\openH}_5$ on $\Pn 4$, $\Hi 1 (F)=V_3$,
$\sCoh 1 F 1=2V_1^{\sharp}$,
$\sCoh 1 F 2=2V_0^{\sharp}$, while the multiplication map
$V^{\ast}\otimes\Hi 1 (F)\to\sCoh 1 F 1$
is given by the projection on the second factor $V^{\ast}\otimes V_3\cong
3V_1\oplus 2V_1^{\sharp}
\to 2V_1^{\sharp}$. One checks easily that
$\xi:=\sum_{i=0}^4 (-1)^i e_i\in\Hi 1 (F)$ has the desired
properties, and thus deduces that there is a Zariski open subset of
elements of $\Hi 1 (F)$ satisfying the two
conditions.\par\medskip

With this choice of $\xi$, the exact sequence $(2.2)$ yields $\hi 0 (E(1))=0$,
$\hi 0 (E(2))=5$, $\hi 1(E(1))=5$ and $\hi 1 (E(m))=0$ for all $m\ge 2$, or
$m\le-1$. Summarizing,
we deduce that $M:=\sCoh 1E{\ast}$ is an artinian module with
Hilbert function $(4,5)$ and with
the desired syzygies:
$$\vbox{
\halign{&\hfil$\,#\,$\hfil\cr
M&\leftarrow&4R&\mapleft{\psi}&15R(-1)&&15R(-2)&&2R(-3)\cr
&&&&&\vbox  to
10pt{\vskip-4pt\hbox{\hss$\nwarrow$}\vss}&\oplus&\longleftarrow&\oplus\cr
&&&&&&12R(-3)&&30R(-4)&\vbox  to
10pt{\vskip-4pt\hbox{\hss$\nwarrow$}\vss}&21R(-5)&\leftarrow&5R(-6)&
\leftarrow&0,\cr
}}$$
where $\psi=\tau\gamma$, with $\tau\in\text{M}_{4,5}(\openC)$
corresponding to the chosen
extension $\xi$, while $\gamma$ is as above.
The two above linear second order syzygies are inherited from
those of $\sCoh 1 F {\ast}$,
and thus they involve as in the Horrocks-Mumford case
two proper Koszul complexes.\par\medskip

As we can check in an example via \cite{7},
the dependency locus of two general sections
in $\sCoh 0E2$ is a smooth surface $S\subset\Pn 4$ with $d=12$, $\pi=13$
and $\chi=3$. The ideal sheaf ${\scr I}_S$ has a resolution of type
$$\exactg {2\scr O} {E(2)}{\id S 5}.\leqno{(2.3)}$$
This description of $S$ is equivalent to that as the
degeneracy locus of a general
morphism $\varphi\in\Hom({\scr E}, Syz_1(\sCoh 1 E {\ast+1}))$.
In particular, $S$ has a minimal
free resolution of type
$$\vbox{%
\halign{&\hfil$\,#\,$\hfil\cr
&&&&3{\scr O}(-5)\cr
0&\leftarrow&{\scr I}_S&\leftarrow&\oplus\cr
&&&&12{\scr O}(-6)&\vbox  to 10pt{\vskip-4pt\hbox{\hss$\nwarrow$}\vss}&30{\scr
O}(-7)&
\longleftarrow&21{\scr O}(-8)&\longleftarrow&5{\scr O}(-9)\leftarrow 0.\cr
}}$$
Dualizing $(2.3)$ we obtain
$$0\longrightarrow\; {\scr
O}(-5)\;\longrightarrow\;E^{\kdual}(-2)\;\longrightarrow\;
2{\scr O}\;\longrightarrow\;\omega_S\;\longrightarrow 0$$
thus $\omega_S$ is globally generated, $p_g=2$, and since the double point
formula yields $K^2=0$
we deduce that $S$ is a minimal proper elliptic surface. Now Le Barz's formula
$(1.2)$ gives
$N_6=10$, hence there are ten $6$-secant lines to $S$ because $E(2)$ is
globally generated
outside a $3$-codimensional set.\qed\enddemo\par

\proclaim {Corollary 2.4} There exist smooth non-minimal $K3$ surfaces
$X\subset\Pn 4$ with
$d=13$, $\pi=16$ which are embedded by a linear system
$$|H|=|H_{min}-7E_0 -\sum_{i=1}^{10} E_i|.$$\endproclaim
\demo{Proof} A minimal proper elliptic surface
$S\subset\Pn 4$ as constructed in Proposition 2.1
can be linked in the complete intersection of
two quintic hypersurfaces to a surface
$X$ with invariants $d=13$, $\pi=16$, $\chi=2$ and a resolution of type
$$\exactg {E^{\kdual}(-2)} {4\scr O} {\id X 5}.\leqno{(2.5)}$$
Smoothness can be checked again in an example. $X$ is cut out
by quintic hypersurfaces, hence there
are no $6$-secant lines this time. On the other hand,
Le Barz's formula gives $N_6=10$, so there exist
$10$ exceptional lines on $X$, namely the $6$-secant lines to $S$.
To describe $X$
we will use adjunction theory. Let $X_1$ denote the
image of $X$ under the adjunction map, defined by $|H+K|$,
and $X_2$ denote the image of $X_1$ under the
adjunction map defined by $|H_1+K_1|$.
We compute the following invariants:\medskip
\settabs \+$S\subset \Pn 4$\qquad&$  H^2=11$\qquad\qquad&$  HK=9$
\qquad\qquad\qquad&$ K^2=-5$\qquad\qquad\qquad&$\pi =11$\cr
\+$X_1\subset\openP^{16}$&$H_1^2=36$&$ H_1K_1=6$&$ K_1^2=-1$&$\pi_1=22$\cr
\+$X_2\subset\openP^{22}$&$H_2^2=47$&$H_2K_2=5$&$ K_2^2=-1+b$&$\pi_2=27,$\cr
\medskip\noindent
where $b$ is the number of $(-1)$-conics on $X$. The Hodge index
theorem gives $K_2^2=-1+b\le 0$, thus
$X$ is either a $K3$, or a proper elliptic surface. Moreover, in case
it is elliptic, $X$ has a $(-1)$ conic or a $(-1)$ cubic and the proper
transform of the canonical
divisor on the minimal model is an elliptic curve of degree $4$ or $5$, while
in case $X$
is a $K3$ surface there is an exceptional rational septic curve on it.
To prove the claim of the corollary it is enough to check that $p_2=1$,
since in the elliptic case Kodaira's formula for the canonical divisor gives
$p_2=\hi 0({\omega_X}^{\otimes 2})\ge 2$.
Dualizing $(2.5)$ we obtain a resolution of $\omega_X$
$$0\longrightarrow\; {\scr O}(-5)\;\longrightarrow\;4{\scr
O}\;\longrightarrow\;E(2)\;
\longrightarrow\;\omega_X\;\longrightarrow 0,$$
and thus one for ${\omega_X}^{\otimes 2}$:
$$0\longrightarrow\; 6{\scr O}\;\longrightarrow\;4E(2)\;\longrightarrow\;
S^2(E)(4)\;\longrightarrow\;{\omega_X}^{\otimes 2}\;\longrightarrow 0.$$
Splitting it up in short exact sequences and using the facts that $\hi
1(E(2))=0$ and
$\hi 0(E(2))=5$ we deduce $p_2=\hi 0(S^2(E)(4))-14$.
Now, from $(2.2)$, $S^2(E)$ can be realized as an extension
$$\exactg {S^2(F)} {S^2(E)} E,$$
while $\hi 0(S^2(F)(4))=0$, $\hi 1(E(4))=0$, and the kernel of the
coboundary morphism $\Hi 0(E(4))\mapright{}\Hi 1(S^2(F)(4))$ has dimension 15,
as one can easily check with $\text {SL}(2,\openZ_5)$-representation theory.
Hence $p_2=1$ and $X$ is a $K3$-surface of the claimed type with a minimal
free resolution
$$
\vbox{%
\halign{&\hfil$\,#\,$\hfil\cr
&&&&4{\scr O}(-5)\cr
0&\leftarrow&{\scr I}_X&\leftarrow&\oplus\cr
&&&&5{\scr O}(-6)&\vbox  to 10pt{\vskip-4pt\hbox{\hss$\nwarrow$}\vss}&16{\scr
O}(-7)&
\leftarrow&10{\scr O}(-8)&\leftarrow&2{\scr O}(-9)\leftarrow 0.\cr
}}\eqno{\qed}$$\enddemo\par

\proclaim {Proposition 2.6} There exist smooth, regular, proper elliptic
surfaces
$S\subset\Pn 4$, with invariants $d=12$, $\pi=14$, $\chi=3$ and embedded by one
of the following linear systems
$$\eqalignno{a)\quad\quad\quad |H|&=|H_{min}-2E_0-\sum_{i=1}^4 E_i|,\cr
b)\quad\quad\quad |H|&=|H_{min}-\sum_{i=1}^5 E_i|.\cr}$$\endproclaim

\demo{Proof} Argueing as in the proof of Proposition 2.1 a
possible Beilinson cohomology table is
$$\qquad\beiliv 2 \ 3 2 \ 1 1 \ $$
Thus we may take this time ${\scr E}:=2{\scr O}(-1)\oplus 3\somega 3 3$ and
${\scr F}:=\ker\psi$,
for some epimorphism $\psi : 2\somega 2 2\oplus \somega 1 1\to {\scr O}$, and
check
for the degeneracy locus of a general morphism $\varphi\in\Hom({\scr E},{\scr
F})$.
Identifying $\Pn 4=\openP(V)$, with $V=\text {span}_\openC(e_0,\dots,e_4)$,
the morphism $\psi$ is induced by a triple $(\psi_{11}, \psi_{12}, \psi_2)$,
where
$\psi_{11}$, $\psi_{12}\in\Dach^2 V$ and $\psi_2\in V$. We check the various
choices for $\psi$. Namely, if \hfill\break
$\alpha)$ $\psi$ is generic, then the associated vector bundle ${\scr F}$ has a
minimal free
resolution of type\par
$$
\vbox{%
\halign{&\hfil$\,#\,$\hfil\cr
0&\leftarrow&{\scr F}&\leftarrow&25{\scr O}(-1)&&&10{\scr O}(-2)&&{\scr
O}(-3)\cr
&&&&&&\vbox  to
10pt{\vskip-4pt\hbox{\hss$\nwarrow$}\vss}&\oplus&\longleftarrow&\oplus\cr
&&&&&&&4{\scr O}(-3)&&4{\scr O}(-4)&\vbox  to
10pt{\vskip-4pt\hbox{\hss$\nwarrow$}\vss}&
{\scr O}(-5)&\leftarrow 0,\cr
}}$$
\smallskip
and a general morphism $\varphi\in\Hom({\scr E},{\scr F})$ gives a smooth
surface
$S_{\alpha}\subset\Pn 4$ with minimal free resolution
$$\vbox{%
\halign{&\hfil$\,#\,$\hfil\cr
0&\leftarrow&{\scr I}_{S_{\alpha}}&\leftarrow&8{\scr O}(-5)&&&7{\scr
O}(-6)&&{\scr O}(-7)\cr
&&&&&&\vbox  to
10pt{\vskip-4pt\hbox{\hss$\nwarrow$}\vss}&\oplus&\longleftarrow&\oplus\cr
&&&&&&&4{\scr O}(-7)&&4{\scr O}(-8)&\vbox  to
10pt{\vskip-4pt\hbox{\hss$\nwarrow$}\vss}&{\scr O}(-9)&
\leftarrow 0\cr}}$$
while, if \hfill\break
$\beta)$ $\psi$ distinguishes a plane, e.g., say $\psi_{11}=0$,
$\psi_{12}=e_0\wedge e_1$ and
$\psi_2=e_2$, then ${\scr F}$ has a resolution
$$
\vbox{%
\halign{&\hfil$\,#\,$\hfil\cr
&&&&25{\scr O}(-1)&&12{\scr O}(-2)&&2{\scr O}(-3)\cr
0&\leftarrow&{\scr
F}&\leftarrow&\oplus&\longleftarrow&\oplus&\longleftarrow&\oplus\cr
&&&&2{\scr O}(-2)&&5{\scr O}(-3)&&4{\scr O}(-4)&\vbox  to
10pt{\vskip-4pt\hbox{\hss$\nwarrow$}\vss}&
{\scr O}(-5)&\leftarrow 0,\cr}}$$
\smallskip
and the generic $\varphi\in\Hom({\scr E},{\scr F})$ degenerates on a smooth
surface $S_{\beta}\subset\Pn 4$
with syzygies
$$\vbox{%
\halign{&\hfil$\,#\,$\hfil\cr
&&&&8{\scr O}(-5)&&9{\scr O}(-6)&&2{\scr O}(-7)\cr
0&\leftarrow&{\scr
I}_{S_{\beta}}&\leftarrow&\oplus&\longleftarrow&\oplus&\longleftarrow&\oplus\cr
&&&&2{\scr O}(-6)&&5{\scr O}(-7)&&4{\scr O}(-8)&
\vbox  to 10pt{\vskip-4pt\hbox{\hss$\nwarrow$}\vss}&{\scr O}(-9)&\leftarrow
0.\cr
}}$$
Smoothness can be checked in examples on a computer via \cite{7}. Also it is
easily seen that all
other choices of $\psi$ lead to singular surfaces, or to
determinantal loci which are not of the expected
codimension. We determine next what
type of surfaces we have constructed.\par\medskip

Let $S_1$ denote the image of $S_{\alpha}$ ($S_{\beta}$ resp.) under the
adjunction map, and let $S_2$ be
the image of $S_1$ under the adjunction map defined by
$|H_1+K_1|$. Then\medskip
\settabs \+$S\subset \Pn 4$\qquad&$  H^2=11$\qquad\qquad&$  HK=9$
\qquad\qquad&$ K^2=-5$\qquad\qquad\quad&$\pi =11$\cr
\+$S_1\subset\openP^{16}$&$H_1^2=35$&$ H_1K_1=9$&$ K_1^2=-5+a$&$\pi_1=23$\cr
\+$S_2\subset\openP^{24}$&$H_2^2=48+a$&$H_2K_2=4+a$&$
K_2^2=-5+a+b$&$\pi_2=27+a,$\cr
\medskip\noindent
where $a$ is the number of $(-1)$-lines and $b$ is the number of $(-1)$-conics
on $S_{\alpha}$
($S_{\beta}$ resp.).\par\smallskip
In case $\alpha)$, the ideal $I_{S_{\alpha}}$ is generated by
quintic hypersurfaces so $S_{\alpha}$ has no $6$-secant
lines. Le Barz's formula $(1.2)$ gives $N_6=4$, hence there are
$4$ exceptional lines, say $E_1, E_2,\dots, E_4$,
on $S_{\alpha}$. Let $S_{\text {min}}$ denote the minimal model of $S_{\alpha}$
and assume first that it is a surface of general
type. Then $S_{\alpha}$ has at least two other exceptional curves $F_1$ and
$F_2$ of degree $\ge$ $2$, and thus there
would exist a curve in $|K_{\text {min}}-F_1|$ of degree $HK_{\text {min}}-2\le
4$ and arithmetic genus $p_a(K_{\text {min}})\ge 2$,
which is a contradiction. It follows that $S_{\alpha}$ is  elliptic, and thus
$K\sim K_{\text {min}}+
\sum_{i=1}^4 E_i + E_0$, where $E_0$ is a $(-1)$ curve of degree $\ge 2$. On
the other side a curve in $|K_{\text {min}}-E_0|$
has arithmetic genus one, so $HE_0\le 3$. We investigate first the case when
$HE_0=3$. Then a curve $D\in|K_{\text {min}}-E_0|$
has degree $4$ and arithmetic genus one,
so it spans only a hyperplane in $\Pn 4$. The residual curve $G\sim H-D$ has
then degree $8$ and genus $9$. We check now that
$G$ lies on a quadric surface in $\Pn 3$. Namely, since Riemann-Roch
gives $\chi(\struct G {2H})=8$ it is enough to show that $\scohos 1 G {2H}\le
1$. This follows from the cohomology of
the exact sequence
$$\eso S {H+D} S {2H} G {2H}$$
since, in the examples above, $\scohos 1 S 1=3$, $\scohos 1 S 2=2$,
$\scohos 2 S {H+D}=\scohos 0 S {K-D-H}=0$,
while the composite multiplication map
$$ \sCohos 1 S H\mapright{D}\sCohos 1 S {H+D}\mapright{H-D}\sCohos 1 S {2H}$$
drops rank at most one on $\Pnd 4$, for a general choice of the morphism
$\varphi\in\Hom({\scr E},{\scr F})$. Now the curve
$D$ lies on two quadrics, thus $\scohid 0 H 4\ge 2$, where $H$ denotes
the hyperplane section of $S_{\alpha}$ cut out
by the $\Pn 3$ of $D$. This is a contradiction, since under the above
assumptions of (minimal) cohomology we have
$\scohid 0 H 4\le \scohid 1 S 3=1$, for all hyperplane sections $H$. Therefore
$E_0$ must be an exceptional conic and
$S_{\alpha}$ a non-minimal elliptic surface embedded by a linear system of type
$a)$, as claimed in the statement
of the proposition.\par\medskip

In case $\beta)$, it is easily seen that the distinguished plane
$\Pi=\openP({\text {span}}_\openC(e_0,e_1,e_2))$ meets $S_{\beta}$
along a plane quintic curve $C$ and the point $P=\openP({\text
{span}}_\openC(e_2))$ outside $C$. Therefore $S_{\beta}$ has infinitely
many $6$-secant lines, namely the pencil of lines in $\Pi$ through $P$, and Le
Barz's formula doesn't apply in this case.
Using the explicit form of the syzygies of ${\scr  I}_{S_{\beta}}$, it is
easily seen that these
are in fact all the $6$-secant lines to $S_{\beta}$. Taking cohomology of the
exact sequences
$$\esi {S_{\beta}} {k-1} {S_{\beta}} k H k \eqno{k=3,4}$$
we observe that $\scohid 1 H 3=3$ for all hyperplane sections $H$ of
$S_{\beta}$, and that $\scohid 1 H 4=1$ if and only
if $P\in H$. Therefore each hyperplane through $P$ contains a plane $\pi$ such
that $\scohid 1 {\pi\cap S_{\beta}} 4=1$.
In particular, the pencil of hyperplanes through $\Pi$ determines a quadric
cone
$$Q=\Bigl\lbrace\;\det\pmatrix
l&m\\
x_3&x_4
\endpmatrix=0\;\Bigr\rbrace,$$
where $l$ and $m$ are suitable linear forms, such that $\scohid 1
{\pi_{(\lambda:\mu)}\cap S_{\beta}} 4=1$
holds for all planes $\pi_{(\lambda:\mu)}=\{\mu l+\lambda m=\mu x_3+\lambda
x_4=0\}$ in one of the rulings of $Q$.
The plane $\Pi$ is obviously a member of the opposite ruling of $Q$, and thus
residual to $C$ in the complete
intersection $S_{\beta}\cap Q$ there is a curve $G\sim 2H-C$ of
degree $19$ and arithmetic genus $23$. On the other side,
the hyperplane sections through $\Pi$ cut, residual to $C$, a pencil $|D|$ with
base point $P$, of curves of
degree $7$ and genus $3$. It follows that the curve $G$ splits as $G=G_1+G_2$,
where $G_1$ is a union of plane curves
contained in planes of the ruling $\pi_{(\lambda:\mu)}$, while $G_2$ is a curve
of degree $14$ which maps down via
projection from the vertex of $Q$ to a complete intersection of type $(7,7)$ on
the quadric surface, which is the base of the cone.
It is easily checked that $G_1$ splits as the union of $5$ exceptional lines,
say $E_1, E_2,\dots,E_5$, on $S_{\beta}$.
Therefore, on the first adjoint surface $S_1$ we obtain $K_1^2\ge 0$ and in
fact, by the Hodge
index theorem, the equality
$K_1^2=0$ holds. We argue further as in case $\alpha)$.
If $S_{\beta}$ were a surface of general type, then
it would contain further an exceptional curve $E$ of degree $\ge 2$ and thus a
curve $N\in|K_{\min}-E|$ would have
degree at most $5$ and arithmetic genus at least $2$. Thus the only case to
exclude is $K_{\min}^2=1$, $HE=2$
and $HN=5$, $p_a(N)=2$. If $N$ spanned only a $\Pn 3$, then the residual curve
$H-N$ would
have degree $7$
and arithmetic genus $8$, which is impossible. Therefore $N$ spans all of $\Pn
4$ and
necessarily splits as $N=A+B$, with $A$ a plane quartic curve and
$B$ a line disjoint from it. But then
$A^2+B^2=N^2=0$, $B^2\le -2$ since $B$ cannot be exceptional,
while the Hodge index theorem yields $A^2\le {16\over 12}$,
which is a contradiction. As claimed, it follows this time that
$S_{\beta}$ is a non-minimal elliptic surface embedded
by a linear system of type $b)$.\qed
\enddemo\par\medskip

\head {3. A family of smooth $K3$ surfaces of degree 14}\endhead
We construct in this chapter an example of a smooth non-minimal $K3$ surface of
degree
$14$ in $\Pn 4$. Currently this is the only known family of smooth non-general
type surfaces
of this degree in $\Pn 4$ . More precisely, we show

\proclaim {Proposition 3.1} There exist smooth non-minimal $K3$ surfaces $S
\subset\Pn 4$, with $d=14$, $\pi=19$, $K^2=-15$, and embedded via
$$|H|=|H_{min}-4E_0-\sum_{i=1}^4 2E_i-\sum_{j=5}^{14} E_j|,$$
where $|H_{min}|$ is very ample on $S_{min}$ of degree $56$ and dimension $29$.
\endproclaim
\demo{Proof} We will discuss two different approaches. First\enddemo
\subhead A syzygy construction\endsubhead
A plausible Beilinson cohomology table for a surface with these invariants is
\bigskip
$$\qquad\beiliw 1 \ 7 7 3 \ \ \ $$
In this case everything is determined by the structure of
$M := \dsum {m\in\openZ}{}\sCohid 2 S {m+4}$.
We will consider the dual module $M^{\ast}$ and  assume that it is
generated, as $R=\openC[x_0,\dots,x_4]$-module, by $\Hom_\openC (\sCohid 2 S 3
, \openC)$.
Thus $M^{\ast}$ has a minimal free resolution of type
$$0\longleftarrow\; M^{\ast}\;\longleftarrow\;
3R(-1)\;\mapleft{\psi}\;8R(-2)\oplus mR(-3)$$
with $m\ge 0$. The morphism $\psi = (\psi_1,\psi_2)$
is given by a $3\times (8+m)$-matrix
with linear entries in $\psi_1$ and quadratic entries in $\psi_2$.
Also, $m>0$ if and only if $\psi_1$ has at least three non-trivial linear
syzygies. However, this doesn't occur for a general choice of $\psi_1$ and
the cokernel of $\psi_1$ is an artinian graded module
with Hilbert function $(3,7,5)$ in
this case. In order to obtain a module with the
desired Hilbert function, it is necessary that the
number of linear syzygies of $\psi_1$ equals $m+2$.
The idea is to start with four planes $P_i=\lbrace l_{i1}=l_{i2}=0\rbrace$,
$i=1,\dots,4$,
and to consider the direct sum of the four Koszul complexes built on
$\lbrace\, l_{i1},l_{i2}\,\rbrace$
$$4R\;\mapleft{\alpha}\; 8R(-1)\;\mapleft{\beta}\; 4R(-2).$$
We take now as $\psi_1$ a
morphism given by three general lines of $\alpha$, i.e.,
$\psi_1 = \gamma\alpha(-1)$, with
$\gamma\in{\text M}_{3,4}(\openC)$ a random matrix, and as $\psi_2$
a general $3\times 2$-matrix with quadratic entries
(since $\psi_1$ has exactly $4$ linear
syzygies). $M^{\ast}:=\coker\psi$ is artinian, with Hilbert function $(3,7,7)$
and with a minimal free resolution of type\par
{\eightpoint
$$ \vbox{
 \halign{&\hfil$\,#\,$\hfil\cr
 0\leftarrow&M^{\ast}\longleftarrow 3R(-1)&&8R(-2)&&4R(-3)\cr
 &&\vbox  to
10pt{\vskip-13pt\hbox{\hss$\overset{\psi}\to\nwarrow$\hss}\vss}
&\oplus&&\oplus\cr
 &&&2R(-3)&\vbox  to 10pt{\vskip-4pt\hbox{\hss$\nwarrow$\hss}\vss}&5R(-4)\cr
 &&&&&\oplus\cr
 &&&&&15R(-5)&\vbox  to
10pt{\vskip-4pt\hbox{\hss$\nwarrow$\hss}\vss}&38R(-6)&\longleftarrow
                              28R(-7)&\longleftarrow 7R(-8)&\leftarrow 0\cr
}}
$$}
We dualize, and set ${\scr F}:={Syz}_2 (M)$. To get a hint for the
second bundle, we compare the syzygies of ${\scr F}$ with Beilinson's spectral
sequence
for ${\scr I}_S$. Namely, the $E_{\infty}$-filtration yields an exact sequence
$$\exactg {{\systruct{-1}}\oplus ({\Hi 0}({\scr F})
\otimes{\scr O})}{\scr F}{\id S 4}.$$
Furthermore $\hi 0 ({\scr F})=15$, so we may take
${\scr E}:={\scr O}(-1)\oplus 15{\scr O}$.
One checks in examples, via \cite{7},
that the degeneration locus of a $\varphi\in
\Hom({\scr E},{\scr F})=\Hom({\systruct{-1}}
\oplus ({\Hi 0}({\scr F})\otimes{\scr O}),{\scr F})=
{\Hi 0}({\scr F}(1))\oplus\Hom({\Hi 0}({\scr F})
\otimes{\scr O},{\scr F})$, given by a general
section and the natural evaluation map,
is a smooth surface $S$ with the desired numerical
invariants and the desired cohomology. The minimal free
resolution of the ideal sheaf of the
surface is of type\par
$$
\vbox{
\halign{&\hfil$\,#\,$\hfil\cr
&&&&4{\scr O}(-5)&&2{\scr O}(-6)\cr
0&\leftarrow&{\scr I}_S&\longleftarrow&\oplus&\longleftarrow&\oplus\cr
&&&&4{\scr O}(-6)&&8{\scr O}(-7)&\vbox  to
10pt{\vskip-4pt\hbox{\hss$\nwarrow$\hss}\vss}
&3{\scr O}(-8)&\leftarrow&0,\cr
}}
$$
and thus the homogeneous ideal is generated by $4$ quintics and $4$ sextics.
Moreover, it follows
from the construction of the module $M^{\ast}$ that the four quintics
containing $S$ intersect
in
$$V((I_S)_{\le 5})=S\cup\bigcup_{i=1}^4 P_i,$$
and a closer look at the syzygies of $M^{\ast}$ shows that $S$ cuts each plane
$P_i$ along
a sextic curve. Hence each of the planes $P_i$ contains an ${\infty}^2$ of
$6$-secant lines,
and in particular Le Barz's formula doesn't apply to this example.\par\medskip
To determine the type of surface we have constructed, one can argue as
follows. All the geometric facts
needed in the sequel will follow from our second construction.
One checks first that
${S\cup\bigcup_{i=1}^4 P_i}$ is an arithmetically
Cohen-Macaulay scheme of degree $18$ and
sectional genus $39$, with syzygies of type\par
$$
\vbox{
\halign{&\hfil$\,#\,$\hfil\cr
0&\leftarrow&{\scr I}_{S\cup\bigcup_{i=1}^4 P_i}&\longleftarrow&4{\scr
O}(-5)&&2{\scr O}(-6)\cr
&&&&&\vbox  to 10pt{\vskip-4pt\hbox{$\nwarrow$}\vss}&\oplus&\leftarrow&0\cr
&&&&&&{\scr O}(-8)\cr
}}
$$
The minors of the above $4\times 2$-submatrix $4{\scr O}(-5)\mapleft{\alpha}
2{\scr O}(-6)$ vanish precisely along an exceptional quartic curve $E_0$ on
$S$. Furthermore,
one may check that there are exactly ten exceptional lines on the surface.
Let now $S_1$ denote the image of $S$ under the
adjunction map, and $S_2$ denote the image of $S_1$ under the map defined by
$|H_1+K_1|$.
Then\medskip
\settabs \+$S\subset \Pn 4$\qquad&$  H^2=11$\qquad\qquad&$  HK=9$
\qquad\qquad\qquad&$ K^2=-5$\qquad\qquad\qquad&$\pi =11$\cr
\+$S_1\subset\openP^{19}$&$H_1^2=43$&$ H_1K_1=7$&$ K_1^2=-5$&$\pi_1=26$\cr
\+$S_2\subset\openP^{26}$&$H_2^2=52$&$H_2K_2=2$&$ K_2^2=-5+b$&$\pi_2=28,$\cr
\medskip\noindent
where $b$ is the number of $(-1)$-conics on $S$. But $K_2^2=-5+b\ge
-H_2K_2=-2$, and there
exists already an exceptional quartic curve $E_0$ on $S$, so the only
possibility is that
$b=4$ and $K_2$ is a $(-1)$ conic on $S_2$. As it turns out, $S=S_{\text
{min}}(p_0,\dots,p_{14})$
is a minimal $K3$ surface blown up in $15$ points and
$$|H|=|H_{\text {min}}-4E_0-\sum_{i=1}^4 2E_i-\sum_{j=5}^{14} E_j|,$$
where $|H_{\text {min}}|$ is a very ample linear system on $S_{\text {min}}$
(the fourth adjoint surface),
defining an embedding $S_{\text {min}}\subset\openP^{29}$, with $\deg S_{\text
{min}}=56$.
\par\bigskip

We want in the sequel to recover an alternative linkage construction for the
above $K3$ surface.
General facts about linkage \cite{33, Prop. 4.1} and \cite{35, Rem. 0.13}
ensure that one can link $(5,5)$ the configuration
$S\cup\bigcup_{i=1}^4 P_i$ to a smooth surface $Y$ of degree $7$ and sectional
genus $6$.
The cohomology of the liaison exact sequence
$$\eso Y {K_Y} {\Sigma_{5,5}} 5 {S\cup\bigcup_{i=1}^4 P_i} 5,$$
where ${\Sigma_{5,5}}$ denotes the complete intersection of the two quintic
hypersurfaces
used in the linkage, gives $p_g(Y)=2$ and $q(Y)=0$ while the double point
formula yields $K_Y^2=0$.
Surfaces with these invariants are classified in \cite{31} and are known to be
arithmetically
Cohen-Macaulay, minimal proper elliptic surfaces. More precisely, $|K_Y|$ is a
pencil without
base points of plane cubic curves; the planes spanned by its members being
those in one
ruling of the determinantal quadric defined by the linear syzygies in
$$\exactg  {2\systruct {-5}}{2{\systruct {-4}}\oplus{\systruct {-2}} } {{\scr
I}_Y}.$$
Once again linkage shows that $Y$ cuts each plane $P_i$ along a conic $C_i$,
which
is necessarily a section of the elliptic fibration, since there are no singular
fibers and
the fibration is by plane curves. In particular, this means that the rank of
the Picard
group of $Y$ is at least $6$, while the Picard number of a generic elliptic
surface of
degree $7$ in $\Pn 4$ is only $2$ by \cite{20}. Therefore $Y$ has to be chosen
carefully
in order to recover $S$ via liaison from the scheme ${Y\cup\bigcup_{i=1}^4
P_i}$.
\par\bigskip

\subhead A linkage construction\endsubhead The above facts suggest
to us the following linkage
construction for this family of $K3$ surfaces. Let $P$, $P_1$, $P_2$, $P_3$,
$P_4$ be five
planes in general position in $\Pn 4$ and denote by $\lbrace p_{ij}\rbrace
=P_i\cap P_j$,
for $1\le i<j\le 4$, the mutual intersection points of the last four of them.

\proclaim {Lemma 3..2} \hfill\break
a) The homogeneous ideal $I_{P\cup\lbrace p_{ij}, {1\le i<j\le 4}\rbrace}$ is
generated by
$3$ quadrics and $4$ cubics.\hfill\break
b) The three quadrics intersect along the plane $P$ and a rational normal
quartic
curve $Q$, which is trisecant to $P$ and goes through the points $p_{ij}$.
\endproclaim

\demo {Proof} The first part follows from the cohomology of the residual
intersection sequences
$$\exactg {\id {\lbrace p_{ij}, i<j\rbrace} {m-1}} {\id {P\cup{\lbrace p_{ij},
i<j\rbrace}} m}
{\struct {\Pn 2} {m-1}}$$
where $m\in\openZ$. For the second part observe that the plane $P$ is linked in
the complete
intersection of two of the hyperquadrics to a rational cubic scroll $T$. If
$H_T\sim C_0+2f$, with
$C_0^2=-1$, $C_0f=1$ and $f^2=0$, is the embedding of the scroll in $\Pn 4$,
then $P\cap T\sim
C_0+f$ is a conic and the third hyperquadric cuts on $T$ the rational normal
quartic curve
$Q\sim C_0+3f$. Now $\sharp P\cap Q =Q(C_0+f)=3$ and the lemma follows.
\qed\enddemo\par

We consider now a general quadric $V\in\sCohid 0 {P\cup\lbrace p_{ij}\rbrace}
2$ and denote with $C_i$
the conics $V\cap P_i$, for $i={1,\dots,4}$.
They intersect pairwise in the points $\lbrace p_{ij}
\rbrace := C_i\cap C_j$, $1\le i<j\le 4$.

\proclaim {Lemma 3.3} There exists a unique rational normal quartic curve $E_0$
which is contained
in $V$, passes through the points $p_{ij}$, $1\le i<j\le 4$, and intersects the
plane $P$ in exactly
one point $p$.\endproclaim

\demo{Proof} The claim is closely related to a theorem of  James \cite{26},
\cite{40}.
Consider the rational map $\gamma :\Pn 4\dashrightarrow\Pn 5$ given by
the quadrics through the rational normal quartic curve $Q$ in lemma 3.2. It is
one to one onto a
smooth hyperquadric $\Omega\subset\Pn 5$, which we identify in the sequel with
the image of
the grassmannian of lines in $\Pn 3$ under the Pl\"ucker embedding. Let
$\widetilde{\openP}^4$ be the
blowing up of $\Pn 4$ along $Q$ and denote by $E$ the exceptional divisor and
by $\tilde\gamma
:\widetilde{\openP}^4\mapright{}\Omega\subset\Pn 5$ the induced morphism. Then
the trisecant planes
of $Q$ are mapped through $\gamma$ to the planes of one generating system, say
$\alpha$-planes, of
the grassmannian $\Omega$, while $E$ is mapped by $\tilde\gamma$ onto a sextic
threefold ruled by
$\beta$-planes. Each of the  $\beta$-planes corresponds to the normal
directions in $\Pn 4$ at points
of $Q$. We remark also that quadric cones through $Q$ are mapped via $\gamma$
to special linear
complexes, i.e., to tangent hyperplane sections of $\Omega$. To fix notations,
let $H\subset\Pn 5$
be the hyperplane corresponding to $V$ and let $\beta_{ij}$, $1\le i<j\le 4$,
be the $\beta$-planes
corresponding to the points $p_{ij}$. Rational normal quartic curves which meet
$Q$ in six points
are represented via $\gamma$ by conics in which $\Omega$ is met by planes.
Thus, in order to prove the lemma,
all we need to check is that there exists exactly one plane in $H$ meeting all
six lines $H\cap\beta_{ij}$
and not contained in the quadric cone $H\cap\Omega$. But this is clear since
the Pl\"ucker embedding of the grassmannian of planes in $\Pn 4$ has degree
$5$, while the planes of
the cone $H\cap\Omega$ describe via the same Pl\"ucker embedding the union of
two conics. The rational
quartic curve $E_0$ represented by this unique plane meets $P$ in one point
because $\gamma$ maps
$P$ to an $\alpha$-plane contained in $H$.\qed\enddemo\par

\proclaim {Lemma 3.4} If $T=P\cup\bigcup_{i=1}^4 C_i$, then its homogeneous
ideal $I_T$ is generated
by $1$ quadric, $2$ cubic and $4$ quartic hypersurfaces.\endproclaim

\demo{Proof} One uses again the residual exact sequences
$$\esi {\cup_{i=1}^4 C_i} {m-1} T m {P\cup(\cup_{i=1}^4 C_i\cap H)} m $$
where $H$ is a general hyperplane through $P$ and $m\in\openZ$, together with
the fact that
$I_{\cup_{i=1}^4 C_i}$ is generated by $1$ quadric and $8$ cubic
hypersurfaces.\qed\enddemo\par

{}From the above lemma it follows that $P$ can be linked in the complete
intersection of $V$ and a
general quartic hypersurface $W\in\sCohid 0 T 4$ to a smooth, minimal proper
elliptic surface
$Y\subset\Pn 4$ with $\deg Y=7$, $\pi(Y)=6$. By construction,
the conics $C_i$ all  all on $Y$.

\proclaim {Lemma 3.5} \hfill\break
a) $C_i^2=-3$ and $K_YC_i=1$ on $Y$, so each conic $C_i$ is a section of the
elliptic
fibration.\hfill\break
b) The planes $P_i$ intersect $Y$ exactly along the conics $C_i$.\endproclaim

\demo{Proof} In any case $K_YC_i\ge 1$ since there are no multiple fibers. On
the other hand, we recall that
the elliptic fibration is cut out on $Y$ by the planes in one of the rulings of
the cone $V$. Thus if
$K_YC_i\ge 2$, then $P_i$ would lie on $V$ and this would contradict our
choices. It follows
that $C_i^2=-3$ and $K_YC_i=1$. Part $b)$ is set theoretically clear by
construction. It is enough
to remark that residual to each conic $C_i$ there is a pencil $|H_Y-C_i|$ of
curves of degree $5$ and
genus $2$, without base points since $(H_Y-C_i)^2=0$.\qed\enddemo\par

We need in the sequel some classical facts of projective geometry.
\proclaim {Proposition 3.6 (Segre)} With any four general planes $P_i$,
$i=1,\dots,4$, there is
associated a uniquely determined fifth plane $P_5$, such that all lines which
meet the first four
planes meet also the fifth.\endproclaim
\demo{Proof} As mentioned above, the Pl\"ucker embedding of the grassmannian of
lines in $\Pn 4$ has
degree $5$, thus the claim follows because the special linear complexes consist
of lines meeting
a given plane. See also \cite{39} or \cite{41}. \qed\enddemo\par

\proclaim {Corollary 3.7 (Segre)\cite{39}} The lines in $\Pn 4$, which meet the
four initial planes,
generate a cubic hypersurface $X$ containing the five planes $P_i$,
$i=1,\dots,5$, and having
singularities (nodes) exactly at the ten points at which the planes meet in
pairs.\endproclaim

\demo{Proof} We briefly recall the arguments in \cite{39}. The first part of
the claim follows from
Bezout's theorem
and from Schubert calculus in $G(1,4)$ since $\sCohid 0 {\cup P_i} 3 =1$, and
since if $l$ is a
line meeting $P_4$ at one point $q$, in which case it is contained in a
hyperplane $H$ through
$P_4$, then there is one line through $q$ which meets $P_1$, $P_2$ and $P_3$,
and there are two
other lines meeting $l$ and the three lines in which $H$ cuts $P_1$, $P_2$ and
$P_3$, respectively.
To check singularities, observe first that residual to a plane $P_i$ in a
general hyperplane
section of $X$ through it, there is a quadric surface containing $4$ skew
lines, thus smooth.
Therefore $X$ has only isolated singularities and an easy argument shows that
these are exactly
the ten points of pairwise intersection of the planes $P_i$.\qed\enddemo\par

A cubic threefold $X\subset\Pn 4$ with the maximum number of ordinary double
points, namely $10$,
is unique up to projective equivalence (cf. \cite{39}, \cite{27}). Its
desingularization $\widetilde X$ is
isomorphic to $\Pn 3$ blown-up in five points $a_1,\dots,a_5$, in general
position. The morphism
$\varphi : {\widetilde X}\to X\subset\Pn 4$ is given by the quadrics through
the five points, while
the nodes are the images of the lines joining any two of the points $a_i$. We
mention in the sequel
some of the properties of this threefold (cf. \cite{39}, \cite{41},
\cite{21}).\par\medskip

The Segre cubic primal $X$ has a symmetrical system of $15$ planes, of which
$5$ correspond to the
exceptional divisors over the points $a_i$, and
$10$ to the planes $P_{ijk}=\varphi({\text {span}}_\openC
(a_i, a_j, a_k))$, for $\lbrace i,j,k\rbrace\subset\lbrace 1,2,3,4,5\rbrace$.
The symmetry of the
the planes resides in the following properties:\hfill\break
- each plane contains four of the nodes, \hfill\break
- each plane is met in lines by $6$ others, namely the plane corresponding to
$a_i$ by the planes
$P_{ijk}$, for all $\lbrace k,j\rbrace\subset\lbrace
1,2,3,4,5\rbrace\setminus\lbrace i\rbrace$,
and the plane $P_{ijk}$ by those corresponding to $a_i$, $a_j$, $a_k$ and
$P_{\alpha,\beta,\gamma}$,
with $\alpha\in\lbrace i,j,k\rbrace$ and $\lbrace\beta,\gamma\rbrace = \lbrace
1,2,3,4,5\rbrace
\setminus\lbrace i,j,k\rbrace$.\par\medskip

We will assume in the sequel that we have chosen the desingularization morphism
$\varphi$ such that
the planes $P_i$, $i=1,\dots,4$, correspond to the exceptional
divisors over $a_i$. Let
now as above $Z=Y\cup\bigcup_{i=1}^4 P_i$. It is a local complete intersection
scheme, outside the
points $P_{ij}$ which are Cohen-Macaulay and where the tangent cone is linked
to a plane
in a complete intersection, and has invariants
$\deg Z=11$, $\pi(Z)=10$, $\chi=3$, $q=0$. We remark here that Hodge index
implies that there
is no smooth surface in $\Pn 4$ with these invariants.
By computing syzygies one shows that $Z$ has a resolution of type
$$\exactg {{2\systruct{-1}}\oplus ({\Hi 0}({\scr G})
\otimes{\scr O})}{\scr G}{\id Z 4}$$
with ${\scr G}:={Syz}_1 (M^\ast)(3)$,
where $M^\ast$ is the graded artinian module in our
first construction. In particular the homogeneous ideal
$I_Z$ is generated by $3$ quintic and $15$
sextic hypersurfaces, and thus we can link $Z$ in the complete intersection of
two quintics
to a surface $S$ with $d=14$, $\pi=19$, $\chi=2$, $q=0$. One checks in examples
via \cite{7}
that $S$ is smooth.

\remark {Remark 3.8} By liaison, each plane $P_i$, $i=1,\dots,4$,
intersects $S$ along
a sextic curve $D_i$, thus each of them contains an $\infty^2$ of $6$-secant
lines to $S$.
\endremark
\proclaim {Lemma 3.9} \hfill\break
a) $E_0$ is an exceptional quartic on $S$.\hfill\break
b) Each of the four planes $P_{ijk}$, with $\lbrace i,j,k\rbrace\subset\lbrace
1,2,3,4\rbrace$, cuts the
surface $S$ along an exceptional conic.\endproclaim

\demo{Proof} The rational normal curve $E_0$ is contained in $V$ and intersects
$W$ in a scheme
of length $16$, of which one point is on $P$. Thus, for general choices, $E_0$
cuts $Z=
Y\cup\bigcup_{i=1}^4 P_i$ along a scheme of length $15+6=21$ and, by Bezout's
theorem, lies on all
quintic hypersurfaces containing $Z$, whence on $S$. We show now that, say
$P_{123}$ cuts $S$
along a conic; the other cases being similar.
Observe first that $P_{123}$ cuts $P_1$, $P_2$
and $P_3$ along the lines pairwise joining the points
$p_{12}$, $p_{13}$, $p_{23}$, while
$P_4$ and $P_5$  both meet this plane at the node $v_{45}$
corresponding to the line through
$a_4$ and $a_5$. For general choices $P_{123}$ meets $Y$
in a scheme of length $7$: $p_{12}$,
$p_{13}$, $p_{23}$ and four extra points.
Let $E_4$ denote the unique conic through these
four points and the node $v_{45}$. It is easily seen that $E_4$ is a
$11$-secant conic to the
configuration $Z$, so by Bezout's theorem it necessarily lies on $S$.\par
By linkage $K_S+((Y\cup\bigcup_{i=1}^4 P_i)\cap S)_S\sim 5H_S$. On the other
hand,
the curve of degree $24$ and arithmetic genus $37$ represented by
$(Y\cap S)_S\sim 5H_Y-K_Y-\sum_{i=1}^4 C_i,$
and the quartic $E_0$ lie both on the quadric cone $V$. It follows that
$K_S+((\bigcup_{i=1}^4 P_i)\cap S)_S\sim 3H_S+E_0$, thus $K_S+\sum_{i=1}^4
D_i\sim (X\cap S)_S
+E_0$. Since $H(K-E_0-\sum_{i=1}^4 E_i)=22-12=10$ and $K^2=-15$ we deduce
easily that $E_i$,
$i={0,\dots,4}$, are exceptional curves on $S$.\qed\enddemo\par

Adjunction and the above lemma show also that $S$ must have $10$ exceptional
lines, thus it is a
non-minimal $K3$ surface embedded by
$$|H|=|H_{min}-4E_0-\sum_{i=1}^4 2E_i-\sum_{j=5}^{14} E_j|.$$

\proclaim {Corollary 3.10} The Segre cubic primal $X$ intersects $S$ along the
union of the
$10$ exceptional lines,  the $4$ exceptional conics and the $4$ plane sextic
curves $D_i$.
\endproclaim

\par\bigskip
A similar linkage construction gives also the following

\proclaim {Proposition 3.11} There exist smooth, non-minimal general type
surfaces $S\subset\Pn 4$
with invariants $d=15$, $\pi=22$, $p_g=3$, $q=0$, $K^2=-6$, and with $9$
exceptional lines.\endproclaim
\demo{Proof} This time one starts with a Castelnuovo surface $Y\subset\Pn 4$,
i.e., with a smooth,
arithmetically Cohen-Macaulay, rational surface with $d=5$, $\pi=2$ and $K^2=1$
(see \cite{8} or \cite{30}).
$Y$ is linked to a plane in the complete intersection of a hyperquadric and a
cubic hypersurface, and
can be represented via the adjunction map as $\openF_1$ blown up in $7$ general
points, thus it is
embedded in $\Pn 4$ by
$$|H_Y|=|4l-2E_0-\sum_{i=1}^7 E_i|.$$
Consider now the following conics on $Y$:
$$\eqalignno{C_0&=3l-2E_0-\sum_{i=1}^6 E_i\cr
C_1&=2l-E_0-E_7-E_1-E_3-E_5\cr
C_2&=2l-E_0-E_7-E_1-E_4-E_6\cr
C_3&=2l-E_0-E_7-E_2-E_4-E_5\cr
C_4&=2l-E_0-E_7-E_2-E_3-E_6.\cr}$$\par
They intersect pairwise in one point and the planes they span, denoted in the
sequel by $P_i$, for
$i\in\lbrace 0,4\rbrace$, intersect the Castelnuovo surface $Y$ only along the
conics $C_i$. The
scheme $Z=Y\cup\bigcup_{i=0}^4 P_i$ is regular, of degree $10$ and sectional
genus $7$, and has a
minimal free resolution of type
$$
\vbox{%
\halign{&\hfil$\,#\,$\hfil\cr
&&&&5{\scr O}(-5)\cr
0&\leftarrow&{\scr I}_Z&\longleftarrow&\oplus\cr
&&&&10{\scr O}(-6)&\vbox  to 10pt{\vskip-4pt\hbox{\hss$\nwarrow$}\vss}&
34{\scr O}(-7)&\longleftarrow&27{\scr O}(-8)&\longleftarrow&7{\scr
O}(-9)\leftarrow 0.\cr}}$$
The five quintics in the ideal intersect along $Z$ and the union of $9$ skew
lines. In fact,
if, according to lemma 3.6 , $Q_i$ denotes the unique Segre cubic hypersurface
containing the planes
$P_{i+1}$, $P_{i+2}$, $P_{i+3}$ and $P_{i+4}$, for $i\in\openZ_5$, then one
checks that $Q_1$, $Q_2$,
$Q_3$ and $Q_4$ each contain $6$ skew lines which are $6$-secant to the
configuration $Z$, while
$Q_0$ contains only $5$ such lines. On the other hand, the five Segre cubics
$Q_i$ cut out an elliptic
quintic scroll $T\subset\Pn 4$ (see \cite{39}, \cite{41, Th.XXXIII, p.278});
each of the planes $P_i$
intersecting it along a cubic curve, section of the ruling. It follows that
there are exactly
$5$ rulings of the scroll which are $6$-secant to $Z$, and thus altogether $9$
skew lines  with
this property. The scheme $Z$ can be linked in the complete intersection of two
quintic hypersurfaces
to a surface $S$ with the desired invariants, having the above $9$ lines as
exceptional curves.
One computes the following cohomology table
\bigskip
$$\qquad\beiliw 3 \ 7 8 4 \ \ \ $$
and a minimal free resolution of type
$$
\vbox{%
\halign{&\hfil$\,#\,$\hfil\cr
&&&&2{\scr O}(-5)\cr
0&\leftarrow&{\scr I}_S&\leftarrow&\oplus\cr
&&&&7{\scr O}(-6)&\vbox  to 10pt{\vskip-4pt\hbox{$\nwarrow$}\vss}&12{\scr
O}(-7)&
\longleftarrow&4{\scr O}(-8)&\leftarrow 0.\cr
}}$$
Finally, we remark that each of the planes $P_i$ intersects $S$ along a sextic
curve, and
thus contains an $\infty^2$ of $6$-secant lines to $S$.\qed\enddemo\par

\remark {Remark 3.12} $S$ and $Z$ are minimal elements in their
even liaison classes.\endremark

\enddocument
\bye